\documentclass[11pt]{amsart}
\usepackage{graphicx}
\usepackage{psfrag}
\usepackage{cite}
\title[Type $D$ WMDs]{Littelmann patterns and Weyl group multiple
Dirichlet series of type $D$} \author{Gautam Chinta}
\address{Department of Mathematics, The City College of CUNY, New
York, NY 10031, USA } \email{chinta@sci.ccny.cuny.edu}

\author{Paul E. Gunnells}
\address{Department of Mathematics and Statistics,
University of Massachusetts,
Amherst, MA 01003, USA
}
\email{gunnells@math.umass.edu}

\date{September 24, 2009}
\thanks{We thank Ben Brubaker and Sol Friedberg for helpful
conversations.  Both authors thank the NSF for support}

\newcommand{\x}{\mathbf{x}}
\newcommand{\g}{\mathfrak{g}}
\newcommand{\p}{\varpi }
\newcommand{\C}{\mathbb{C}}
\newcommand{\Q}{\mathbb{Q}}
\newcommand{\psia}{\ensuremath{{\it \Psi}}}
\renewcommand{\O}{\mathcal{O}}
\newcommand{\ba}{\overline{a}}
\newcommand{\m}{\mathbf{m}}
\newcommand{\s}{\mathbf{s}}
\renewcommand{\c}{\mathbf{c}}

\DeclareMathOperator{\SO}{SO}

\DeclareMathOperator{\SL}{SL}
\DeclareMathOperator{\GL}{GL}
\DeclareMathOperator{\fSL}{\mathfrak{sl}}
\theoremstyle{plain}
\newtheorem{conjecture}{Conjecture}
\theoremstyle{definition}

\newtheorem{example}{Example}
\newtheorem{definition}{Definition}

\begin{document}
\begin{abstract}
We formulate a conjecture for the local parts of Weyl group multiple
Dirichlet series attached to root systems of type $D$.  Our conjecture
is analogous to the description of the local parts of type $A$ series
given by Brubaker, Bump, Friedberg, and Hoffstein \cite{wmd3} in terms of
Gelfand--Tsetlin patterns.  Our conjecture is given in terms of
patterns for irreducible representations of even orthogonal Lie
algebras developed by Littelmann \cite{l}.
\end{abstract}

\maketitle

\section{Introduction}
We begin with some notation.  Let $\Phi$ be a reduced
root system of rank $r$ and $n$ a positive integer.  Let $F$ be a
number field containing the $2n$-th roots of unity.  Let $S$ be a set
of places of $F$ containing the archimedean places and those that
ramify over $\Q$, as well as sufficiently many more places to ensure
that the ring of $S$-integers $\O_{S}$ has class number $1$.  Let $\m
= (m_{1},\dotsc ,m_{r})$ be a fixed nonzero tuple of elements of
$\O_{S}$.  Let $\s = (s_{1},\dotsc ,s_{r})$ be an $r$-tuple of complex
variables.

Given the data above, one can form a \emph{Weyl group multiple
Dirichlet series}.  This is a Dirichlet series in the $r$ variables
$s_{i}$ with a group of functional equations isomorphic to the Weyl
group $W$ of $\Phi$.  More precisely, one can define a set of
functions of the form
\[
Z (\s ; \m , \psia) = Z_{\Phi}^{n} (\s; \m , \psia) = \sum_{\c}
\frac{H (\c; \m )\psia (\c)}{\prod |c_{i}|^{s_{i}}},
\]
where each $c_{i}$ ranges over nonzero elements of $\O_{S}$ modulo
units, $\psia$ is taken from a certain finite-dimensional complex
vector space $\Omega$ of functions on $(F_{S}^{\times})^{r}$, and $H$
is an important function we shall say more about shortly.  Then the
collection of all such $Z$ as $\psia$ ranges over a basis of $\Omega$
satisfies a group of functional equations isomorphic to $W$ with an
appropriate scattering matrix.  For more about why Weyl group multiple
Dirichlet series are interesting objects, as well as a discussion
about the basic framework for their construction, we refer to
\cite{wmd1}.

The heart of the construction of $Z$ is the function $H$.  This
function must be carefully defined to ensure that $Z$ satisfies the
correct group of functional equations.  The heuristic of \cite{wmd1}
dictates how to define $H$ on the \emph{powerfree} tuples $\c, \m $
(those tuples such that the product $c_{1}\dotsb c_{r}m_{1}\dotsb
m_{r}$ is squarefree).  Moreover, it is further specified in
\cite{wmd1} how the values of $H$ on the \emph{prime power} tuples $\c
= (\p^{k_{1}},\dotsc ,\p^{k_{r}})$, $\m = (\p^{l_{1}},\dotsc
,\p^{l_{r}})$, where $\p \in \O_{S}$ is a prime, determine $H$ on all
tuples.

Thus, writing $\ell$ for a tuple of nonnegative integers
$(l_{1},\dotsc ,l_{r})$ and letting $\p^{\ell}$ denote the tuple
$(\p^{l_{1}},\dotsc ,\p^{l_{r}})$, the construction of $Z$ reduces to
understanding the multivariate generating function
\begin{equation}\label{eqn:genfun}
N (x_{1},\dotsc ,x_{r}; \ell ) := \sum_{k_{i}\geq 0} H
(\p^{k_{1}},\dotsc ,\p^{k_{r}}; \p^{\ell }) x_{1}^{k_{1}}\dotsb
x_{r}^{k_{r}}.
\end{equation}

At present there are two different approaches to understanding the
generating function \eqref{eqn:genfun}, and thus to constructing 
Weyl group multiple Dirichlet series.  Both are related to characters
of representations of the semisimple complex Lie algebra attached to
$\Phi$.  Let $\omega_{i}, i=1,\dotsc ,r$ be the fundamental weights of
$\Phi$ and let $\theta$ be the
strictly dominant weight $\sum (l_{i}+{1})\omega_{i}$.

\begin{itemize}
\item The \emph{Gelfand--Tsetlin} approach \cite{wmd3, wmd4, wmd5a},
which works for $\Phi = A_{r}$, gives formulas for the coefficients
$H(\p^{k_{1}},\dotsc ,\p^{k_{r}}; \p^{\ell })$.  These formulas are
written in terms of Gauss sums and statistics extracted from
Gelfand--Tsetlin patterns for the representation of $\fSL_{r+1} (\C)$
of lowest weight $-\theta $.
\item The \emph{averaging} approach \cite{qmds, cfg, wmdsn, nmds},
which works for all $\Phi$, uses a ``metaplectic'' deformation of the
Weyl character formula to construct a rational function with known
denominator, whose numerator is then taken to define $N$.
\end{itemize}

Both approaches have their advantages and limitations.  The
Gelfand--Tsetlin construction gives very explicit formulas for $H$,
formulas that (remarkably) are uniform in $n$ and that lead to a
direct connection with the global Fourier coefficients of Borel
Eisenstein series on the $n$-fold cover of $\SL_{r+1}$ \cite{wmd5},
but suffers from the obvious disadvantage that it only works for type
$A$.  The averaging approach, on the other hand, works for all $\Phi$,
quickly leads to the definition of $Z$, yet has the drawback that it
seems difficult to get similarly explicit formulas for the
coefficients of $N$.  By combining recent work of Chinta--Offen
\cite{co} and McNamara \cite{mcnamara}, we know that in type $A$ the
two definitions of $N$ coincide, although it seems difficult to give a
direct combinatorial proof.

This note arose from our attempts to understand the Gelfand--Tsetlin
approach to \eqref{eqn:genfun}.  In the course of studying
\cite{wmd3}, it became plain to us that the most suitable language to
understand the constructions in \cite{wmd3} is that of Kashiwara's
\emph{crystal graphs}, as encoded in the generalization of the
Gelfand--Tsetlin basis due to Littelmann \cite{l}, which we call
\emph{Littelmann patterns}.  Indeed, the definitions in \cite{wmd3}
become much more transparent when phrased in terms of these patterns.

To test the relevance of this observation, we decided to try to
formulate a Littelmann analogue of the Gelfand--Tsetlin construction
when $\Phi$ is a root system of type $D$.  The main result of this
note is thus Conjecture \ref{conj:ppartconj}, which explicitly
describes the generating function $N (x_{1},\dotsc ,x_{r} ; \ell )$
for the $\p$-part of the type $D$ Weyl group multiple Dirichlet series
constructed using the averaging method.  We remark that for $n=1$,
Conjecture \ref{conj:ppartconj} gives a type $D$ analogue of a theorem
of Tokuyama \cite{tokuyama}.

We have some limited evidence for the truth of Conjecture \ref{conj:ppartconj}.

First, for $D_{2}\simeq A_{1}\times A_{1}$, the conjecture is easily
seen to be true.

Next, we have tested the conjecture for $D_{3}$ when $n\leq 4$ and for
$D_{4}$ when $n\leq 2$, by computing the $\p$-parts by averaging for
many tuples $\ell$ and comparing with the predictions of Conjecture
\ref{conj:ppartconj}.  In all cases there was complete agreement.
Note that $D_{3}\simeq A_{3}$, so the $\p$-part of the $D_{3}$-series
has already been described explicitly using the results of
\cite{wmd3}, and in this guise has already been compared extensively
with $\p$-parts constructed by averaging.  Nevertheless, agreement in
rank $3$ between $\p$-parts constructed using Conjecture
\ref{conj:ppartconj} and using averaging is a nontrivial check, since
$D_{3}$ Littelmann patterns are quite different from $A_{3}$
Gelfand--Tsetlin patterns.

Finally, recently Brubaker and Friedberg have computed the global
Whittaker coefficients of Eisenstein series on covers of $\GL_{4}$ by
inducing from the parabolic subgroup of type $\GL_{2}\times \GL_{2}$
\cite{bf}.  Their computations---which build on earlier work of
Bump--Hoffstein \cite{bh} and are the first attempts to extend the
results of \cite{wmd5} beyond type $A$ and to work with other
parabolic subgroups---express the Whittaker coefficients in terms of
certain exponential sums.  In the course of their work Brubaker and
Friedberg found that the integrals can be broken up in accordance with
the decomposition of $H (\p^{k_{1}},\dotsc ,\p^{k_{r}}; \p^{\ell })$
given by Conjecture \ref{conj:ppartconj}, and that if one does so the
contributions to the global Whittaker coefficient exactly agrees with
Conjecture \ref{conj:ppartconj}.  We find this connection between
Eisenstein series and $\p$-parts to be strongly convincing evidence of
the correctness of Conjecture \ref{conj:ppartconj}.

\section{Littelmann patterns}
Let $\g$ be the simple complex Lie algebra of type $D_{r}$, in other
words the Lie algebra of the group $\SO_{2r} (\C)$.  Let $\theta$ be a
dominant weight of $\g$ and let $V_{\theta}$ be the irreducible
$\g$-module of highest weight $\theta$.  In \cite[\S 7]{l} Littelmann
describes a way to index a basis of $V_{\theta}$ using patterns that
are analogous to the classical Gelfand--Tsetlin patterns for the Lie
algebra of $\SL_{r} (\C)$.  In this section we recall his
construction.

\begin{figure}[htb]
\psfrag{1}{$1$}
\psfrag{2}{$2$}
\psfrag{3}{$3$}
\psfrag{4}{$4$}
\psfrag{5}{$5$}
\psfrag{6}{$6$}
\begin{center}
\includegraphics[scale=0.4]{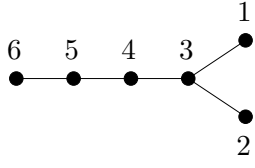}
\end{center}
\caption{The diagram for $D_{6}$\label{fig:diagram}}
\end{figure}

First we label vertices of the Dynkin diagram of $\g$ with the
integers $1,\dotsc ,r$.  We label the upper node of the right prong $1$,
the lower node of the prong $2$, the node at the elbow of the prong
$3$, and then the remaining nodes increase from $4$ to $r$, reading
right to left (Figure \ref{fig:diagram}).  We remark that this is
not the standard labelling by Bourbaki, which begins with $1$ at the
left of the diagram.

A pattern $T$ for $D_{r}$ consists of a collection of integers
$a_{i,j}$, where $1\leq i\leq r-1$ and $i\leq j\leq 2r-2$.  We picture
$T$ by drawing the integers placed in $r-1$ rows of centered boxes.
The first row contains $2r-2$ boxes, the second $2r-4$ boxes, and so
on down to the $(r-1)$-st row, which contains $2$ boxes.  The integers
are placed in the boxes so that $a_{i,i}$ is placed in the leftmost
box of the $i$th row, and then the remaining integers $a_{i,j}$ are
put in the boxes in order as $j$ increases.  We define an involution
on each row by $\ba_{i,j} = a_{i,2r-1-j}$.

To index a weight vector in $V_{\theta}$, there are two sets of
inequalities the $a_{i,j}$ must satisfy.  The first is independent of
$\theta$: in each row we must have
\begin{equation}\label{eqn:admissible}
a_{i,i}\geq a_{i,i+1}\geq \dotsb \geq a_{i,r-2} \geq a_{i,r-1},
a_{i,r} \geq a_{i,r+1}\geq \dotsb \geq a_{i, 2r-1-i}\geq 0,
\end{equation}
or, using the bar notation, 
\[
a_{i,i}\geq a_{i,i+1}\geq \dotsb \geq a_{i,r-2} \geq a_{i,r-1},
\ba_{i,r-1} \geq \ba_{i,r-2}\geq \dotsb \geq \ba_{i, i}\geq 0
\]
In other words, the $a_{i,j}$ are weakly decreasing in the rows, with
the exception that no comparison is made between $a_{i,r-1}$ and
$a_{i,r}$.  Both of these entries, however, are required to be $\leq
a_{i,r-2}$ and $\geq a_{i,r+1}$.

\begin{definition}
A pattern $T$ is \emph{admissible} if $T$ satisfies
\eqref{eqn:admissible} for all $i$.  
\end{definition}

\begin{figure}[htb]
\psfrag{1}{$1$}
\psfrag{2}{$2$}
\psfrag{3}{$3$}
\psfrag{4}{$4$}
\psfrag{5}{$5$}
\psfrag{6}{$6$}
\begin{center}
\includegraphics[scale=0.45]{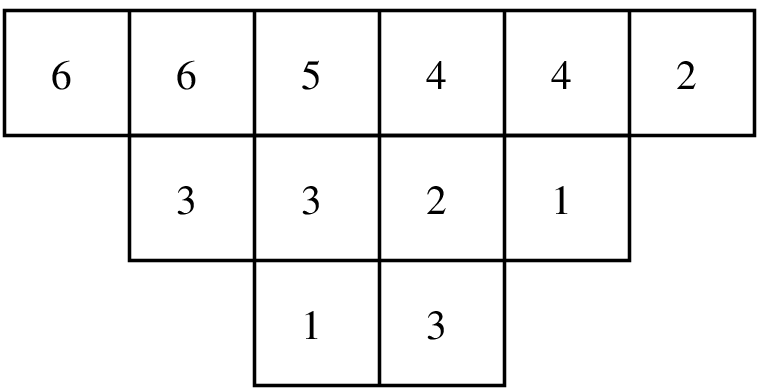}
\end{center}
\caption{\label{fig:example}}
\end{figure}

Figure \ref{fig:example}
shows an admissible pattern for $D_{4}$.  

The next set of inequalities involves the highest weight $\theta$.
Write 
\[
\theta = \sum m_{k}\omega _{k},
\]
where the $\omega_{k}$ are the fundamental weights.  Then an
admissible $T$ will
correspond to a weight vector in $V_{\theta}$ if $T$ satisfies
\begin{align}
\label{eq1}\ba_{i,j}&\leq m_{r-j+1} + s (\ba_{i,j-1}) - 2s (a_{i-1,j})+s
(a_{i-1,j+1})\quad \text{for $j\leq r-2$},\\
a_{i,j}&\leq  m_{r-j+1} + s (a_{i,j+1}) - 2s (\ba_{i,j})+s
(\ba_{i,j-1})\quad \text{for $j\leq r-2$},\\
a_{i,r-1} &\leq m_{2} + s (\ba_{i,r-2})-2t (a_{i-1,r-1}), \text{and}\\
\label{eqlast}a_{i,r} &\leq m_{1} + s (\ba_{i,r-2})-2t (a_{i-1,r}),
\end{align}
where we write for $j<r-1$
\begin{align*}
s (\ba_{i,j}) &= \ba_{i,j} + \sum_{k=1}^{i-1} (a_{k,j}+\ba_{k,j}),\\
s (a_{i,j}) &= \sum_{k=1}^{i} (a_{k,j}+\ba_{k,j}),\\
s (a_{i,r-1}) &=s (\ba_{i,r-1}) = \sum_{k=1}^{i} a_{k,r-1} + a_{k,r},\\
t (a_{i,r-1}) &= \sum_{k=1}^{i}a_{k,r-1}, \quad t (a_{i,r}) = \sum_{k=1}^{i}a_{k,r}.
\end{align*}

\begin{definition}
A pattern $T$ is \emph{$\theta$-admissible} if $T$ is admissible and its
entries satisfy \eqref{eq1}--\eqref{eqlast}.
\end{definition}

Note that the inequalities for the $i$th row only involve the entries
of $T$ on the $i$th and $(i-1)$st rows.  Moreover when ordered in
terms of increasing $i$, there is a unique inequality in which a given
entry $a_{i,j}$ appears on the left.

\begin{definition}\label{def:critical}
We say that an entry in a $\theta$-admissible pattern is
\emph{critical} if this first inequality is actually an equality.
\end{definition}

To complete our discussion of Littelmann patterns, we must assign a
weight $\lambda (T) $ to each pattern $T$.  This is a vector $\lambda
(T) = (\lambda_{1},\dotsc ,\lambda_{r})$ of nonnegative integers,
where 
\[
\lambda_{k} = \begin{cases}
\sum_{i=1}^{r-1}(a_{i,r+1-k}+\ba_{i,r+1-k})&k=3,4,\dots ,r\\
\sum_{i=1}^{r-1}a_{i,r-2+k}&k=1,2.
\end{cases}
\]
We write $|\lambda | = \lambda_{1}+\dotsb +\lambda_{r}$.
In our conjecture, if a pattern $T$ occurs for the twist $\theta =
\sum m_{i}\omega_{i}$, it will contribute to the coefficient of
$\x^{\lambda (T)} := x_{1}^{\lambda_{1}}\dotsb x_{r}^{\lambda_{r}}$ in
the numerator $N (\x , \ell )$, where $\ell = (l_{1},\dotsc ,l_{r})$
and $l_{i}=m_{i}-1$.  For instance, the pattern in Figure
\ref{fig:example} contributes to the coefficient of
$x_{1}^{9}x_{2}^{9}x_{3}^{14}x_{4}^{8}$, with $x_{1}$ corresponding to
left middle column of three entries and $x_{2}$ to the right middle
column of three entries.

\section{The decorated graph of a pattern} 

Let $T$ be a $\theta$-admissible pattern.  We want to associate to $T$
a graph $\Gamma (T)$.  The graph $\Gamma (T)$ will also potentially be
endowed with \emph{decorations}, which will be circled vertices.  The
vertices of $\Gamma (T)$ correspond to the entries of $T$; the graph
will have at least one connected component for each row of $T$.

\begin{figure}[htb]
\psfrag{1}{$a_{1,1}$}
\psfrag{2}{$a_{1,2}$}
\psfrag{3}{$a_{1,3}$}
\psfrag{4}{$a_{1,4}$}
\psfrag{5}{$a_{1,5}$}
\psfrag{6}{$a_{1,6}$}
\psfrag{7}{$a_{1,7}$}
\psfrag{8}{$a_{1,8}$}
\psfrag{9}{$a_{1,9}$}
\psfrag{10}{$a_{1,10}$}
\begin{center}
\includegraphics[scale=0.4]{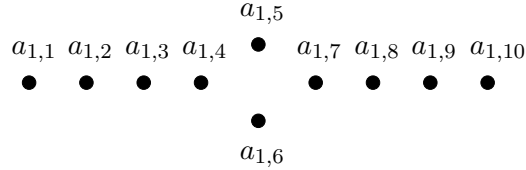}
\end{center}
\caption{The vertices for the top row of $D_{6}$\label{fig:vertices}}
\end{figure}

We begin by describing how each row determines a subgraph.  Consider
the $i$th row of $T$. Each entry $a_{i,j}$ in this row gives an
vertex.  We draw the corresponding vertices in a row, with the two
vertices in the middle corresponding to the incomparable entries
$a_{i,r-1}, a_{i,r}$ entries arranged vertically.  For definiteness we
assign $a_{i,r-1}$ to the top vertex and $a_{i,r}$ to the bottom
vertex.  See Figure \ref{fig:vertices} for the arrangement for the top
row of a pattern for $D_{6}$.

Now join two vertices by an edge if they appear consecutively in the
inequalities \eqref{eqn:admissible}, are equal, and are comparable in
\eqref{eqn:admissible}.  Note that we do not join the vertices
corresponding to $a_{i,r-1}, a_{i,r}$ by an edge if they happen to be
equal, since they are not comparable in \eqref{eqn:admissible}.  This
gives a graph for this row.  We then do the same for each row of $T$.
The result is $\Gamma (T)$ without decorations.  

Certain symmetric connected components that arise in the construction
of $\Gamma (T)$ will play a special role in our conjecture:

\begin{definition}\label{def:ml}
Let $T$ be an admissible pattern and suppose $a_{i,j} = \ba_{i,j}$ for
some $i,j$ with $j\not =r-1,r$.  Then the component of
$\Gamma (T)$ containing $a_{i,j}, \ba_{i,j}$ is called a
\emph{multiple leaner}.  If in addition $a_{i,j-1}\not =
a_{i,j}$ and $\ba_{i,j-1}\not =\ba_{i,j}$ then we say the multiple
leaner is \emph{symmetric}. We define the \emph{length} $l (C)$ of a
symmetric multiple leaner to be half the number of its vertices.
\end{definition}

\begin{figure}[htb]
\begin{center}
\includegraphics[scale=0.4]{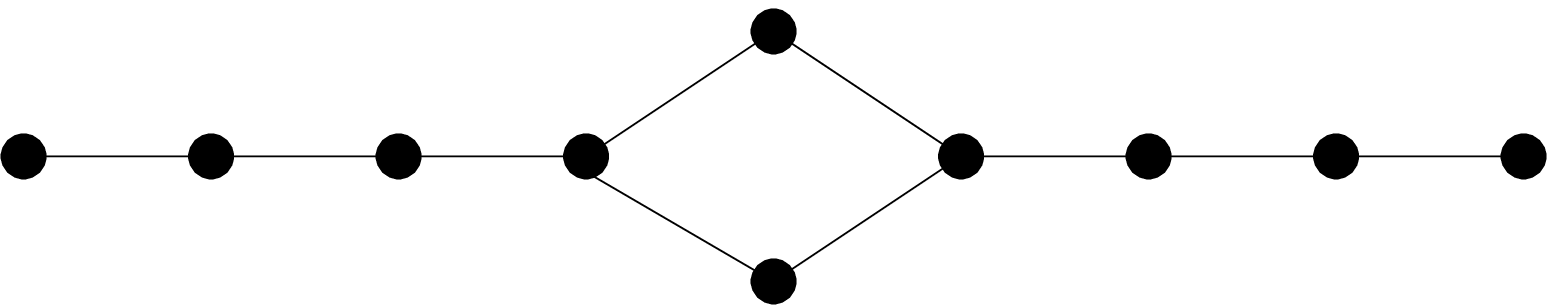}
\end{center}
\caption{\label{fig:multiple}}
\end{figure}

The term \emph{leaning} is inspired by \cite{wmd3}; see also
\S\ref{s:lsc}.  Figure \ref{fig:multiple} shows an example of a
symmetric multiple leaner of length $5$, when all the entries in the
top row of a pattern for $D_{6}$ are equal.  Note that the minimal
length of a symmetric multiple leaner is $2$, and that multiple
leaners can appear in patterns for $D_{3}$, but not for $D_{2}$.

To complete the construction of $\Gamma (T)$ we must describe how to
add the decorations.  This is very simple: we circle each vertex whose
corresponding entry is critical in the sense of Definition
\ref{def:critical}.

\begin{figure}[htb]
\psfrag{1}{$1$}
\psfrag{2}{$2$}
\psfrag{3}{$3$}
\psfrag{4}{$4$}
\psfrag{5}{$5$}
\psfrag{6}{$6$}
\begin{center}
\includegraphics[scale=0.45]{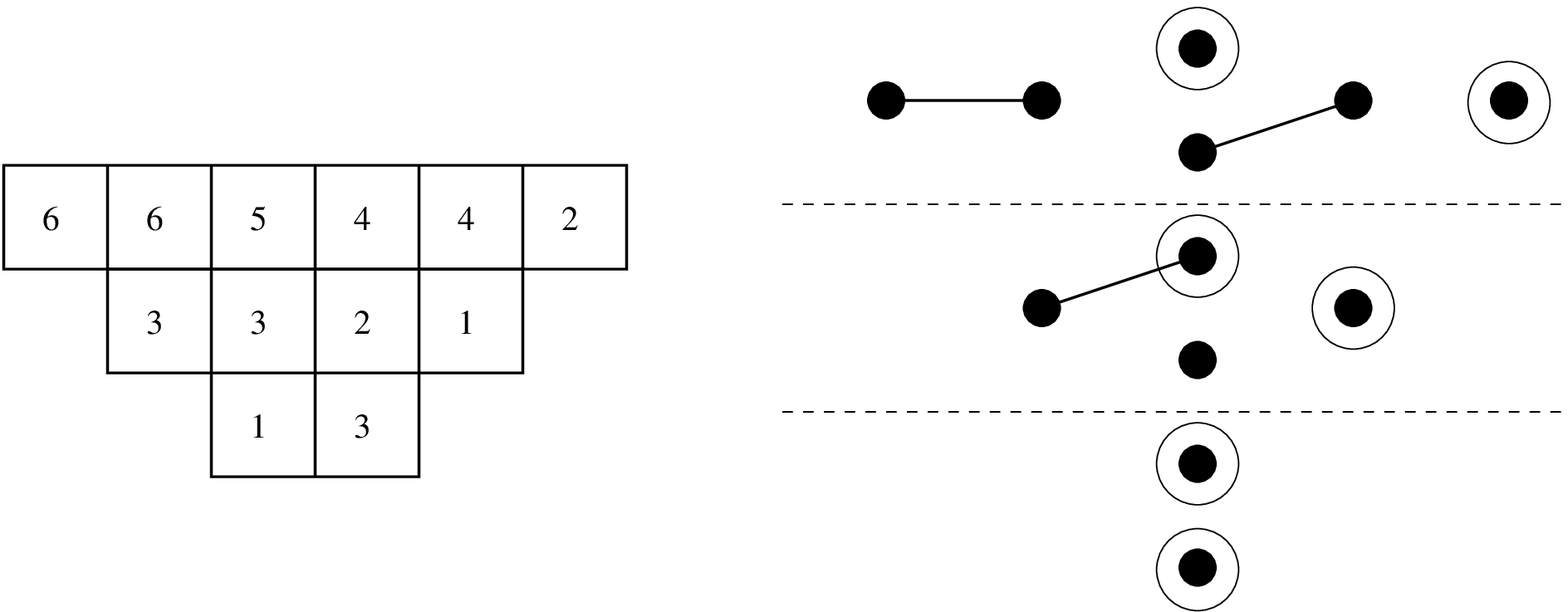}
\end{center}
\caption{\label{fig:decorated}}
\end{figure}

Figure \ref{fig:decorated} shows an example of building the decorated
graph of the Littelmann pattern in Figure \ref{fig:example}.  We
assume that a highest weight $\theta $ has been specified so that the
circled vertices in the graph correspond to critical entries.

\section{Strictness}

In the following for $k=1,\dotsc ,n-1$ we write $g_{k}$ for the Gauss
sum $g (\p ^{k-1},\p ^{k})$ (see for example \cite{wmdsn} for the
definition of the Gauss sums).  For convenience we extend the notation
and define $g_{0}$ to be $-1$.  It is also convenient to define
$g_{m}$ for $m\geq n$ by $g_{m} = g_{k}$, where $m = k$ mod $n$ and
$k=0,\dots ,n-1$.  We let $p$ be the norm of $\p$.

\begin{figure}[htb]
\begin{center}
\includegraphics[scale=0.6]{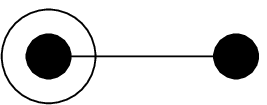}
\end{center}
\caption{\label{fig:nonstrict}}
\end{figure}

In \cite{wmd3} certain patterns for a given weight are discarded and
do not contribute to the relevant coefficient of $N$; such patterns
are called \emph{strict} in \cite{wmd3}.  In type $A$ strictness
corresponds to an easily stated property for Gelfand--Tsetlin
patterns.  If one interprets the definition of strictness in
\cite{wmd3} in terms of type $A$ Littelmann patterns, one sees that a
type $A$ pattern is nonstrict exactly when
\begin{itemize}
\item an entry is simultaneously $0$ and critical, or
\item there are two  adjacent entries that are equal, with the left
entry critical.  
\end{itemize}
We take these to be our definition for type $D$ patterns as well:

\begin{definition}
A type $D$ Littelmann pattern $T$ is called \emph{strict} if the
following conditions hold:
\begin{itemize}
\item No component of $\Gamma (T)$ contains a vertex with a circled $0$.
\item No component of $\Gamma (T)$ that is not a multiple-leaner
contains a subgraph of the form shown in Figure \ref{fig:nonstrict}
(in this figure, the rightmost vertex is less than the left vertex in
the partial order from \eqref{eqn:admissible}).
\end{itemize}
\end{definition}

Note that the subgraph from Figure \ref{fig:nonstrict} is allowed to
appear in multiple-leaners.

\section{Leaning and standard contributions} \label{s:lsc} Let $T$ be
a strict pattern, and let $\Gamma (T)$ be the associated decorated
graph.  For any connected component $C$ of $\Gamma (T)$, let $y_{C}$
be the \emph{rightmost} vertex, in the sense of the order induced by
the inequalities \eqref{eqn:admissible}.  If $C$ has two rightmost
vertices, meaning that it is in the $i$-th row and contains entries
$a_{i,r-2} = a_{i,r-1} = a_{i,r} \not = a_{i,r+1}$, then we define the
rightmost vertex to be the vertex corresponding to $a_{i,r-1}$, that
is, the upper vertex in Figure \ref{fig:vertices}.

\begin{definition}
Let $T$ be a pattern and $\Gamma = \Gamma (T)$ the associated
decorated graph.  Fix $n$ and let $y$ be an entry of $T$.  We define
the \emph{standard contribution} $\sigma (y)$ by the following rule:
\begin{itemize}
\item If the vertex corresponding to $y\not =0$ is uncircled, then we put
$\sigma (y) = 1-1/p$ if $n$ divides $y$ and $\sigma (y) = 0$
otherwise.
\item If the vertex corresponding to $y\not =0$ is circled, then we put
$\sigma (y) = g_{k}/p$, where $y = k \bmod n$ and $k=0,\dotsc ,n-1$.
\end{itemize}
Note that $\sigma (y)$ depends on $n$ and $\theta$, even though we
omit them from the notation. 
\end{definition}

We are almost ready to state our conjecture.  There is one more
phenomenon that plays a role, namely \emph{leaning}.  Essentially,
leaning means that if entries are consecutive and equal in a
Littelmann pattern $T$, where consecutive means adjacent in
\eqref{eqn:admissible}, then only one should contribute to the
corresponding coefficient of $N (\x ; \ell )$.  This is why we
introduce the graph $\Gamma (T) $.  Its connected components keep
track of these equalities among entries.

Thus we are led to consider contributions of the connected components
of $\Gamma (T)$, not just the entries.  There is further slight twist that
the contribution of a multiple leaning component is different from
that of all other components:

\begin{definition}
Let $C$ be a connected component of $\Gamma (T)$.  The \emph{standard
contribution $\sigma (C)$ of $C$} is defined as follows:
\begin{itemize}
\item If $C$ is not a multiple leaner, then we put $\sigma (C) =
\sigma (y_{C})$, where $y_{C}$ is the
rightmost entry of $C$.
\item If $C$ is a multiple leaner that is not symmetric, let $y_{C}$
be the entry on the endpoint of its shorter leg.  Then we define
$\sigma (C) = \sigma (y_{C})$.
\item If $C$ is a symmetric multiple leaner, then let $y_{C}$ be its rightmost
entry $a_{i,j}$ and $\upsilon_{C}$ (upsilon =
Greek $y$) to be the entry $a_{i,j-1}$.  Then we define
\[
\sigma (C) = \begin{cases}
\sigma (y_{C}) (1-1/p^{l (C)})&\text{if $y_{C}$ is uncircled,}\\
\sigma (y_{C})\sigma (\upsilon_{C}) (1/p^{l (C)-1})&\text{if $y_{C}$
is circled,}
\end{cases}
\]
where $l (C)$ is defined to the half the number of vertices of $C$
(Definition \ref{def:ml}).
\end{itemize}
\end{definition}

We are now ready to state our conjecture:

\begin{conjecture}\label{conj:ppartconj}
Let $N (\x ; \ell) = \sum_{\lambda} a_{\lambda}\x^{\lambda}$ be the
$\p$-part constructed by averaging \cite{qmds, wmdsn} for the Weyl
group multiple Dirichlet series $Z_{\Phi}^{n} (\s; \m , \psia)$.  Then
we have
\begin{equation}\label{eqn:form}
a_{\lambda} = p^{|\lambda|}\sum_{T} \prod_{C \subset \Gamma (T)} \sigma ({C}), 
\end{equation}
where the sum is taken over all strict patterns $T$ of weight
$\lambda$ and with highest weight $\theta = \sum (l_{i}+1)\omega_{i}$,
and the product is taken over the connected components of $\Gamma
(T)$.
\end{conjecture}

\begin{example}
Suppose the pattern in Figure \ref{fig:example} appears for a highest
weight $\theta$ such that the decorated graph appears in Figure
\ref{fig:decorated}.  Suppose $n=2$.  Then the contribution of this
pattern to the coefficient of $x_{1}^{9}x_{2}^{9}x_{3}^{14}x_{4}^{8}$
will be
\[
p^{40} \Bigl(1-\frac{1}{p}\Bigr)^{3} \Bigl(-\frac{1}{p}\Bigr)
\Bigl(\frac{g_{1}}{p}\Bigr)^{5}.
\]
\end{example}

\begin{example}
We consider another example for $n=2$.  Suppose the twisting parameter
is $\ell = (0,1,2,0)$, which corresponds to the highest weight $\theta
= \omega_{1}+2\omega_{2}+3\omega_{3}+\omega_{4}$.  We will compute the
coefficient $a_{\lambda}$ of the monomial
$\x^{\lambda }=x_{1}^{10}x_{2}^{10}x_{3}^{17}x_{4}^{10}$.  Note that
$|\lambda | = 47$.  

There are $27$ Littelmann patterns that we must consider.  Six of
these patterns are nonstrict, for instance the pattern shown in Figure
\ref{fig:ns}.  Of the remaining $21$, only $2$ give nonzero
contributions; these patterns $T_{1}$, $T_{2}$ appear in Figures
\ref{fig:p1}--\ref{fig:p2}.  Note that Figure \ref{fig:p2} contains a
multiple leaner of length $2$.  All of the other $19$ patterns have an
odd entry that is not circled, and thus have a connected component in
$\Gamma (T)$ with standard contribution equal to zero.

Each vertex in $\Gamma (T_{1})$ is its own connected component.  We
see $3$ uncircled even nonzero entries, $5$ circled even nonzero
entries, and $3$ circled odd entries.  Thus $T_{1} $ contributes
\[
p^{47} \Bigl(1-\frac{1}{p}\Bigr)^{3} \Bigl(-\frac{1}{p}\Bigr)^{5}
\Bigl(\frac{g_{1}}{p}\Bigr)^{3}
\]
to $a_{\lambda}$.

The pattern $T_{2}$ has a multiple leaner $C$ of length $l (C) = 2$.
Its rightmost entry $y_{C}$ is circled, and the entry $\upsilon_{C}$
is uncircled.  We have $\sigma (y_{C}) = -1/p$, $\sigma (\upsilon_{C})
= (1-1/p)$.  These appear in \eqref{eqn:form} multiplied by the
additional factor $1/p$ to account for the length of $C$.  Each of the
remaining vertices is its own connected component, and we have no
uncircled nonzero evens, $4$ circled nonzero evens, and $3$ circled
odds.  Thus $T_{2}$ contributes
\[
p^{47} \Bigl(-\frac{1}{p}\Bigr)^{4}\Bigl(\frac{g_{1}}{p}\Bigr)^{3}\Bigl(-\frac{1}{p}\Bigr)\Bigl(1-\frac{1}{p}\Bigr) \Bigl(\frac{1}{p}\Bigr)
\]
to $a_{\lambda}$.  After simplifying we find 
\[
a_{\lambda} = -p^{36} \left(p^3-2 p^2+2 p-1\right)g_{1}^3,
\]
in agreement with the $\p$-part from \cite{qmds}.
\begin{figure}[htb]
\psfrag{0}{$0$}
\psfrag{1}{$1$}
\psfrag{2}{$2$}
\psfrag{3}{$3$}
\psfrag{4}{$4$}
\psfrag{5}{$5$}
\psfrag{6}{$6$}
\psfrag{7}{$7$}
\psfrag{8}{$8$}
\psfrag{9}{$9$}
\psfrag{10}{$10$}
\begin{center}
\includegraphics[scale=0.45]{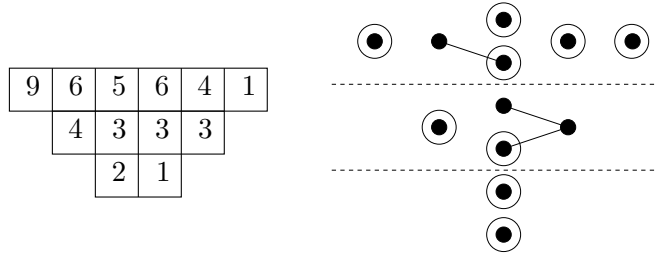}
\end{center}
\caption{A nonstrict pattern\label{fig:ns}}
\end{figure}

\begin{figure}[htb]
\psfrag{0}{$0$}
\psfrag{1}{$1$}
\psfrag{2}{$2$}
\psfrag{3}{$3$}
\psfrag{4}{$4$}
\psfrag{5}{$5$}
\psfrag{6}{$6$}
\psfrag{7}{$7$}
\psfrag{8}{$8$}
\psfrag{9}{$9$}
\psfrag{10}{$10$}
\begin{center}
\includegraphics[scale=0.45]{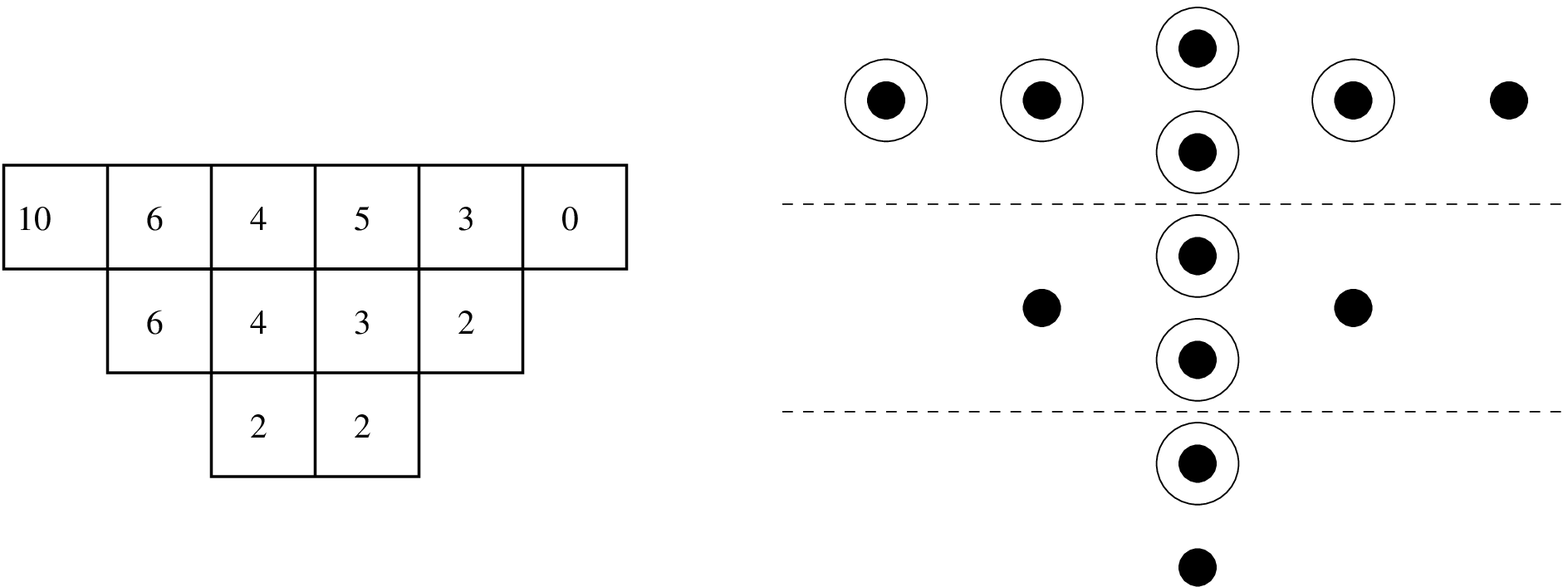}
\end{center}
\caption{$T_{1}$\label{fig:p1}}
\end{figure}

\begin{figure}[htb]
\psfrag{0}{$0$}
\psfrag{1}{$1$}
\psfrag{2}{$2$}
\psfrag{3}{$3$}
\psfrag{4}{$4$}
\psfrag{5}{$5$}
\psfrag{6}{$6$}
\psfrag{7}{$7$}
\psfrag{8}{$8$}
\psfrag{9}{$9$}
\psfrag{10}{$10$}
\begin{center}
\includegraphics[scale=0.45]{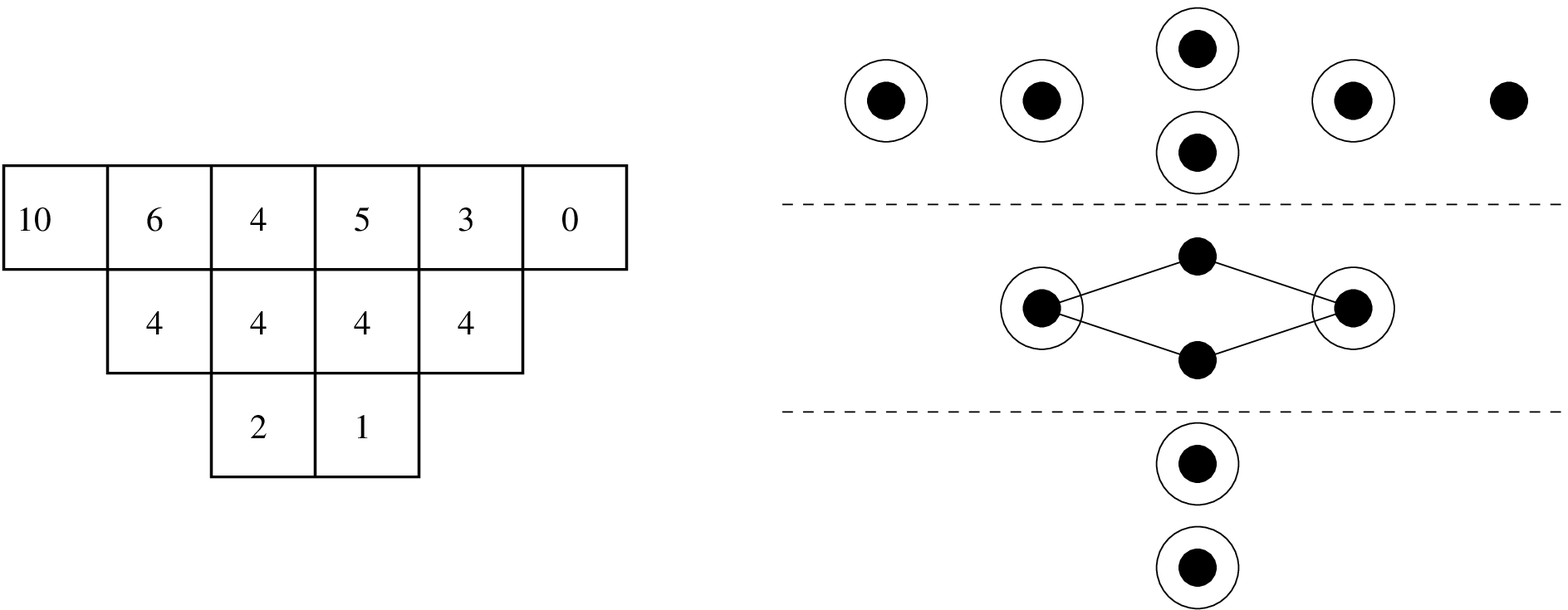}
\end{center}
\caption{$T_{2}$\label{fig:p2}}
\end{figure}

\end{example}

\begin{example}
We conclude by describing an example for $D_{4}$ when $n=1$.  We put
$\ell = (0,0,0,0)$ (the ``untwisted'' case), so that the highest
weight is $\omega_{1}+\omega_{2}+\omega_{3}+\omega_{4}$.  The
polynomial $N (\x ;\ell)$ is supported on $601$ monomials.  There are
$4096$ Littelmann patterns to consider, $2216$ of which are nonstrict.
The remaining patterns each give a nonzero contribution to $N$.  The
resulting polynomial can be written succinctly as
\[
N (\x ; \ell) = \prod_{\alpha > 0} (1-p^{d (\alpha)-1}\x^{\alpha}),
\]
where the product is taken over the positive roots.  Here $d (\alpha)
= k_{1}+k_{2}+k_{3}+k_{4}$ if $\alpha$ is the linear combination of
simple roots $k_{1}\alpha_{1}+k_{2}\alpha_{2}+k_{3}\alpha_{3}+k_{4}\alpha_{4}$, and
$\x^{\alpha}$ refers to the monomial $x_{1}^{k_{1}}x_{2}^{k_{2}}x_{3}^{k_{3}}x_{4}^{k_{4}}$.
\end{example}

\bibliographystyle{amsplain_initials}
\bibliography{d4conjecture}
\end{document}